\documentclass[12pt,a4paper]{article}

\usepackage{amsmath,amssymb,amsthm}
\usepackage[mathcal]{euscript}

\usepackage[dvips]{graphicx}
\usepackage[dvips]{color}

\newtheorem{theorem}{Theorem}
\newtheorem{lemma}[theorem]{Lemma}

\newtheorem{corollary}[theorem]{Corollary}

\newtheorem*{remark}{Remark}

\newcommand{\GL}{G_\linel}

\newcommand{\Galph}{G_\alpha}

\newcommand{\linel}{\mathfrak{L}}
\newcommand{\rootpartition}{\mathfrak{P}}
\newcommand{\pointP}{\mathcal{P}}
\newcommand{\hatt}{\hat {\ } }
\newcommand{\spaceS}{\mathcal{S}}

\begin{document}

\pagestyle{plain}
\pagenumbering{arabic}

\title{Linear spaces with significant characteristic prime}
\author{Nick Gill}

\maketitle

\begin{abstract}
Let $G$ be a group with socle a simple group of Lie type defined over the finite field with $q$ elements where $q$ is a power of the prime $p$. Suppose that $G$ acts transitively upon the lines of a linear space $\spaceS$. We show that if $p$ is {\it significant} then $G$ acts flag-transitively on $\spaceS$ and all examples are known.

\end{abstract}

\section{Background and statement of result}

A {\it linear space} $\spaceS$ is an incidence structure of points and lines such that any two points are incident with exactly one line. Also $\spaceS$ is {\it non-trivial} provided any point is incident with at least two lines and any line is incident with at least two points; all linear spaces considered in this paper will be presumed to be non-trivial. A {\it flag} is a pair $(\alpha, L)$ where $\alpha$ is a point incident with a line $L$.

Let $\spaceS$ be a finite linear space admitting an automorphism group $G$ which is transitive on lines. Then $\spaceS$ is said to have parameters $b$ (the number of lines), $v$ (the number of points), $k$ (the number of points incident with a line) and $r$ (the number of lines incident with a point). 

Camina, Neumann and Praeger \cite{cnp} have defined a prime $p$ to be {\it significant} for the space $\spaceS$ if it divides into $(b,v-1)$. They then show that if $P$ is a Sylow $p$-subgroup of $G$ and $\Galph$ is a point-stabilizer in $G$ then $\Galph\geq N_G(P)$ \cite[Lemma 6.1]{cnp}.

The finite linear spaces which admit a flag-transitive almost simple group have been classified in \cite{kleidman4, saxl}. As part of the program to extend this classification to those linear spaces which admit a line-transitive almost simple group we prove the following theorem:

\begin{theorem}\label{theorem:main}
Suppose that a group $G$ has socle a group of Lie type of characteristic $p$. Suppose furthermore that $G$ acts transitively upon the lines of a linear space $\spaceS$ with significant prime $p$. Then $G$ acts transitively upon the flags of $\spaceS$ and we have one of the following examples:
\begin{itemize}
\item $U_3(q)\leq G \leq P\Gamma U(3,q)$ and $\spaceS$ is a Hermitian unital.

\item ${^2G_2(q)}\leq G \leq Aut({^2G_2(q)})$ and $\spaceS$ is a Ree unital.

\end{itemize}
\end{theorem}

The remainder of this paper will be occupied with a proof of Theorem \ref{theorem:main}. The suppositions given in Theorem \ref{theorem:main} will be assumed from here on.

\section{A reduction to simplicity}

Observe that, by \cite[Lemma 6.1]{cnp} mentioned above, a point-stabilizer $\Galph$ must contain a parabolic subgroup of the socle of $G$. We can use this fact along with the notion of {\it exceptionality} to immediately simplify our task. 

Let $G_0$ be a normal subgroup in a group $G$ which acts upon a set $\pointP$. Then $(G,G_0,\pointP)$ is called {\it exceptional} if the only common orbital of $G_0$ and $G$ in their action upon $\pointP$ is the diagonal (see \cite{gms}). Then we have the following result:

\begin{lemma}\cite[Lemma 26]{gill2}
Suppose a group $G$ acts line-transitively on a finite linear space $\spaceS$; suppose furthermore that $G_0$ is a normal subgroup which is not line-transitive on $\spaceS$; finally suppose that $|G:G_0|=t$, a prime.

Then either $\spaceS$ is a projective plane or $(G,G_0,\pointP)$ is exceptional where $\pointP$ is the set of points in $\spaceS$.
\end{lemma}

Now consider a pair $(G,\spaceS)$ satisfying the suppositions of Theorem \ref{theorem:main}. Then $\spaceS$ is not a projective plane since the finite projective planes are precisely the finite linear spaces with no significant prime. Thus if $G$ contains a normal subgroup $G_0$ of index a prime $t$ which is not line-transitive on $\spaceS$ then $(G,G_0,\pointP)$ is exceptional.

However all of the exceptional triples of this form are enumerated in \cite[Theorem 1.5]{gms}. In all cases a point-stabilizer does not contain a parabolic subgroup of the socle of $G$. We can conclude from this that our socle itself is transitive on the lines of $\spaceS$.

In fact, referring to \cite{csk}, we see that if the socle of $G$ has Lie rank 1 then it acts $2$-transitively upon its parabolic subgroups. Thus the socle of $G$ is $2$-transitive upon the points of $\spaceS$ and hence is transitive on the flags of $\spaceS$ (c.f. \cite{bdd}). Then, by \cite{saxl}, the actions listed in Theorem \ref{theorem:main} are the only examples.

Thus for the remainder of this paper we add the following suppositions to those mentioned in Theorem \ref{theorem:main}: 
\begin{itemize}
\item	We suppose that $G$ is simple;

\item	We suppose that $G$ has Lie rank greater than 1.
\end{itemize}
We will show that these suppositions lead to a contradiction. We will do this by taking $\Galph$ to be a parabolic subgroup of $G$ and then examining potentional line stabilizers, $\GL$.

\subsection{Group theory notation}

In our  use of the theory of groups of Lie type we will use the notation of Carter\cite {carter}. For $G$ a Chevalley group we have the standard subgroups $B, U, H, N$ and the associated Weyl group $W$. We write $\Phi$ and $\Pi$ be the set of roots, and the set of fundamental roots respectively, associated with $G$. 

For $G$ a twisted simple group, consider $G$ as a subgroup of $G^*$ the untwisted simple group. Let $\Phi$ and $\Pi$ for the set of roots, and the set of fundamental roots respectively, associated with $G^*$ and take $\rho$ to be the non-trivial symmetry of the Dynkin diagram. Take $W^1$ to be the Weyl group of $G$, a subgroup of $W$, the Weyl group of $G^*$. The subgroups $U^1, V^1, H^1$ and $N^1$ are defined as usual. Write $\rootpartition$ for the partition of $\Pi$ into $\rho$-orbits.

We will sometimes precede the structure of a subgroup of a projective group with $\hatt$ which means that we are giving the structure of the pre-image in the corresponding universal group. An integer $n$ denotes a cyclic group of order $n$, while $[n]$ denotes an arbitrary soluble group of order $n$.

\section{The point stabilizer is non-maximal}\label{section:nonmaximal}

\begin{lemma}
Suppose that $G$ is a simple Chevalley group acting on a linear space $\spaceS$ with $\Galph$ a non-maximal parabolic subgroup of $G$. Then $\GL$ is a parabolic subgroup of $G$ and $p$ is not significant.
\end{lemma}
\begin{proof}
Let $\Phi^+$ be the set of positive roots associated with $G$ so that
$$U=\prod_{s\in \Phi^+} X_s.$$
For $r\in\Pi$ be a fundamental root define the group 
$U_r=\prod_{s\in \Phi^+\backslash\{r\}} X_s.$

Now suppose that $\Galph$ is the parabolic subgroup $P_J$ where $J$ is a subset of $\Pi$ the set of fundamental roots. Since $\Galph$ is non-maximal in $G$ we know that at least two fundamental roots, say $s$ and $t,$ do not lie in $J$.

For $s$ a fundamental root recall the standard homomorphism $\phi_s$ from $SL(2,q)$ into $\langle X_s, X_{-s} \rangle$. Then 
$$n_s:=\phi_r\big(
\begin{array}{cc}
0 & 1 \\
-1 & 0 
\end{array}
\big)$$
Now $n_s$ is an involution lying outside of $\Galph$ but which normalizes $U_s$ inside of $\Galph$. Hence $U_s$ fixes at least two points and hence the line between them. So $\GL$ contains a $G$-conjugate of $U_s$. Similarly $\GL$ contains a $G$-conjugate of $U_t$. In fact $\GL$ contains a $G$-conjugate of $U_s:H$ and $U_t:H$.

Now consider a Sylow $p$-subgroup of $\GL$. For some choice of $\linel$ this lies inside $U$. Now observe that, since $G=BNB$ and since both $U_s$ and $U$ are normal in $B$,
\begin{eqnarray*}
&&U_s^g< U \\
&\implies& b_1nb_2U_sb_2^{-1}n^{-1}b_1^{-1}<U \ {\rm where} \ g=b_1nb_2\\
&\implies& nU_sn^{-1}<U \\
&\implies& nU_sn^{-1}=U_s.
\end{eqnarray*}

Thus $U$ only contains one $G$-conjugate of $U_s$ and one $G$-conjugate of $U_t$, namely themselves. Furthermore they generate $U$. Thus $\GL$ contains $B=U:H$ as required.

Now $p$ does not divide into $b$ and so $p$ is not significant.
\end{proof}

\begin{lemma}
Suppose that $G$ is a twisted simple group acting on a linear space $\spaceS$ with $\Galph$ a non-maximal parabolic subgroup of $G$. Then $\GL$ is a parabolic subgroup of $G$ and $p$ is not significant.
\end{lemma}
\begin{proof}
Let $J$ be a $\rho$-orbit of $\Pi$. Then observe that
$$U_J^1=\big(\prod_{r\in \Phi^+, r\not\in\Phi_J^+} X_r\big) \cap U^1$$
is a subgroup of $U_1$ which is normalized by $\langle X^1_{\Phi_J^+}, X^1_{\Phi_J^-}\rangle$.

Let $w^1_J$ be the element in $W^1$ which maps every positive root of $\Phi_J$ to a negative root of $\Phi_J$. Then, by \cite[Proposition 13.5.2]{carter}, there exists $n^1_J\in N^1$ which maps onto $w^1_J$ in the natural way. Now $w^1_J$ can be thought of as a reflection and $(n_J^1)^2\in H^1$.

Now suppose that $\Galph$ lies inside the parabolic subgroup $P_{\rootpartition\backslash\{J,K\} }$ where $J$ and $K$ are distinct $\rho$-orbits of $\Pi$. Then $n^1_J$ and $n^1_K$ do not lie in $\Galph$. By the same argument as above this means that $\GL$ contains a $G$-conjugate of $U_J^1:H^1$ and $U_K^1:H^1$.  

As before consider a Sylow $p$-subgroup of $\GL$. For some choice of $\linel$ this lies inside $U^1$. Furthermore just as before $U^1$ only contains one $G$-conjugate of $U^1_J$ and one $G$-conjugate of $U^1_K$ and these generate $U$. Thus $\GL$ contains $B^1=U^1:H^1$ and we have a contradiction.
\end{proof}

\section{The point-stabilizer is maximal}\label{section:maximal}

In this section take $G$ to be a Chevalley group. Our argument generally translates in a straightforward way to the twisted groups and so we will not repeat it; we will comment on any deviations as we proceed. For convenience we note that, by trivial combinatorial arguments, $G={^2F_4(2)'}$ cannot act line-transitively upon our linear space $\spaceS$.

Take $r\in\Pi$ and suppose that $\Galph = P_J$ where $\Pi = J\cup \{r\}$. By the argument in the previous section  
%
%
it is clear that $\GL\geq U_rL_{\Pi\backslash K}$ where $L_{\Pi\backslash K}$ is the Levi complement of the parabolic group $P_{\Pi\backslash K}$ and $K=\{r\}\cup K'$ where
$$K'=\{ {\rm fundamental \ roots \ which \ are \ not \ orthogonal \ to \ } r\}.$$ 


Observe first of all that, for the Chevalley groups, if $\GL$ contains any $p$-element $h$ from 
$$\langle U_r, X_r, X_{-r}\rangle \backslash U_r$$ 
then $\GL\geq \langle h, U_rH\rangle = B^g$ for some $g\in G$. This is a contradiction.

For the twisted groups this argument does not work in all cases. We need to show that $U^1_J:H$ is maximal in all conjugates of the Borel of which it is a subgroup. It is sufficient to show that $H$ acts transitively upon set of the non-identity elements of $X_J^1$. We refer to \cite[Tables 2.4 and 2.4.7]{gorenstein} to see that this is only true when $X_J^1$ is of type I, II, III and VI as listed there. The cases we have excluded are when $G={^2A_n}(q)$, $n$ even, with $\Galph=\hatt [q^{\frac{n^2+4n}{4}}]:GL_\frac{n}{2}(q^2)$; and when $G={^2F_4}(q)$ with $\Galph=[q^{22}]:GL_2(q^2)$.

Now we will investigate the possibility that there exists $g\in \GL \backslash (P_{\Pi\backslash K'}\cap\GL)$. Suppose that this is the case. Since we have a $BN$ pair we can write $g=u_1 n_w u$ where $u_1,u\in U$ and $n_w\in N$ maps onto $w\in W$ under the natural epimorphism. In fact, since $\GL\geq U_rH$ we can assume that $g= x_r(t)n_wx_r(u)$ where $t,u$ are elements of the finite field of order $q$. 

Now suppose that $w(r)\neq \pm r$ (and note that then $w^{-1}(r)\neq \pm r$). We seek to prove the following

\begin{equation}
g^{-1}U_rg \ \cap \langle U_r, X_r, X_{-r}\rangle \not\leq U_r.\label{equation:overlap}
\end{equation}

Clearly we can replace $g$ by $n_w$ since $X_r$ normalizes $U_r$ and $\langle X_r, X_{-r}\rangle$. So we are required to prove

$$n_w^{-1}U_rn_w \ \cap \langle U_r, X_r, X_{-r}\rangle\not\leq U_r.$$

Since $w(r)\neq\pm r$ we know that, for some $s\in\{r,-r\}$,
$$n_wX_sn_w^{-1}<U_r.$$

This implies (\ref{equation:overlap}) and so there exists a $p$-element in $\GL$ lying in 
$$\langle U_r, X_r, X_{-r}\rangle\backslash U_r.$$
This element will normalize $U_r$ and so $\GL\geq B$. This is a contradiction.

Thus if there exists $g\in \GL \backslash (P_{\Pi\backslash K'}\cap\GL)$ then we can take $g=u_1 n_w u$ as before and $w(r)=\pm r$. In fact, just as before, we can without loss of generality assume that $g=x_r(t)n_wx_r(u)$.

Now suppose that for all $s$, adjacent fundamental roots of $r,$ we have $w(s)$ in $\Phi^+\cup \Phi_{\Pi\backslash K}^-$. Since $\GL> L_{\Pi\backslash K}$ we can assume that $w(s)$ is positive for all fundamental roots not equal to $r$. But then, by \cite[Theorem 2.2.2]{carter}, $w=w_r$ or $w=1$ (see also \cite[Lemma 13.1.3]{carter} for the twisted case). However $\GL$ also contains $n_r$ and so we can assume that $g=x_{\pm r(t)}x_r(u)$. In this case though $g\in P_{\Pi\backslash K'}$ which is a contradiction.

Thus there exists $s$ an adjacent fundamental root of $r$ such that $w(s)$ is negative. Define $h:= gx_s(v)g^{-1}$. As before we can suppose that $h=x_r(v_1)n_{w_1}x_r(v_2).$

Now observe that $g\in \langle X_r, X_{-r}\rangle N_N(\langle X_r, X_{-r}\rangle )$. Suppose that $h$ also lies in $\langle X_r, X_{-r}\rangle N_N(\langle X_r, X_{-r}\rangle )$. Then this would imply that 
$$x_s(v)\in \langle X_r, X_{-r}\rangle N_N(\langle X_r, X_{-r}\rangle ).$$
This is clearly impossible, see \cite[Corollary 8.4.4, Proposition 13.5.3]{carter}.

Thus $h\not\in \langle X_r, X_{-r}\rangle N_N(\langle X_r, X_{-r}\rangle ).$  This implies that $w_1(r)\neq\pm r$. Furthermore since $w(s)\not\in\Phi^+\cup \Phi_{\Pi\backslash K}^-$, $h\not\in P_{\Pi\backslash K'}$. Then we can apply the same argument to $h$ as we applied to $g$ above. This will lead us to conclude that $\GL\geq B$ which is a contradiction.

This leads to the following result:
\begin{lemma}\label{lemma:parabolicadj}
Suppose that $G$ is a Chevalley group with $\Galph=P_{\Pi\backslash r}$. Then
$$U_rL_{\Pi\backslash K}\leq \GL\leq P_{\Pi\backslash K'}.$$

Suppose alternatively that $G$ is a twisted group with $\Galph=P_{\rootpartition\backslash J}$. Suppose further that $G\neq {^2F_4}(q)'$ and $G\neq {^2A_n(q)}$, $n$ even. Then
$$U^1_JL_{\rootpartition\backslash K}\leq \GL\leq P_{\rootpartition\backslash K'}$$
where $K=J\cup K'$ and $K'$ is the set of orbits of fundamental roots in $\rootpartition$ which contain roots not orthogonal to some root in $J$.
\end{lemma}

We record the following lemma of Saxl:

\begin{lemma}\cite[Lemma 2.6]{saxl}\label{lemma:saxl}
If $X$ is a group of Lie type of characteristic $p$ acting on cosets of a maximal parabolic subgroup then there is a unique subdegree which is a power of $p$ except where $X$ is one of $PSL_n(q)$, $P\Omega^+_{2m}(q)$ ($m$ odd) or $E_6(q)$.
\end{lemma}

For the moment let us exclude the exceptions listed in these two lemmas; then Lemma \ref{lemma:saxl} suggests that if $\Galph=P_r$ then $\GL$ contains some $G$-conjugate of $L_r$. This clearly contradicts Lemma \ref{lemma:parabolicadj}. Note also that even in the listed exceptions of Lemma \ref{lemma:saxl} many of the maximal parabolic subgroups have a unique subdegree which is a power of $p$.

\subsection{The twisted exceptions}

We consider the exceptional cases listed in Lemma \ref{lemma:parabolicadj}. In fact we need only consider when $(G,\Galph)$ is one of (${^2A_n}(q)$, $\hatt [q^{\frac{n^2+4n}{4}}].GL_\frac{n}{2}(q^2))$, $n$ even; or $({^2F_4}(q)$, $[q^{22}]:GL_2(q^2)), q^2=2^{1+2a}, a\geq 1.$

In both cases Lemma \ref{lemma:saxl} still applies. Furthermore if $\Galph=P_{\rootpartition\backslash J}$ then $\GL\geq U_J^1$ and so $b$ divides $|X_J|v$.

Consider the unitary case. Write $\Galph=P_{\rootpartition\backslash\{b\}}$ where $b$ is the missing root class. Now $|\GL|$ is divisible by $\frac{P_{\rootpartition\backslash\{b\}}}{q^3}$ and we examine the maximal subgroups of ${^2A_n}(q)$ (\cite{kl}) to find that, unless $(n,q)\in\{(9,2),(11,2)\}$, $\GL<P_{\rootpartition\backslash\{b\}}$ for some $\linel$. The exceptions can be eliminated by trivial counting arguments.

By the work in Section \ref{section:nonmaximal},
$$U_b^g< U\implies U_b^g=U_b.$$
Thus if we choose $\linel$ such that there exists $P\in Syl_p\GL$ with $P<U$ then $U_b<\GL<U.L_{\rootpartition\backslash\{b\}}$. Now $\GL$ contains a Levi complement of $P_{\rootpartition\backslash\{b\}}$ so, in particular, contains an element $g:=un_a$. Here $u\in U$ and $n_a$ is an element of $N$ which when mapped to the Weyl group is the reflection in root class $a$ where $a$ is adjacent to $b$. Without loss of generality we can assume that $g=x_b(t)n_a$. Then
$$gX_b g^{-1}= x_b(t)n_aX_b n_a^{-1}x_b(t)^{-1}=x_b(t)X_{w_a(b)}x_b(t)^{-1}<U_b.$$
Since $U_b<\GL$ this implies that $X_b<\GL$ which is a contradiction.

Now when $G={^2F_4}(q)$ it is clear that $|\GL|$ is divisible by $\frac{|P_J|}{q^4}$. Number the root classes of $G$ as corners of a 16-gon. Then the fundamental root classes are 1 and 8; $\Galph=P_{\{8\}}$. Examining the subgroups of $G$ (\cite{malle}), $\GL<P_{\{8\}}$ for some $\linel$.

By the work in Section \ref{section:nonmaximal}
$$U_1^g< U\implies U_1^g=U_1.$$
Thus if we choose $\linel$ such that there exists $P\in Syl_p\GL$ with $P<U$ then $U_1<\GL<U.L_{\{8\}}$. Now $\GL$ contains a Levi complement of $P_{\{8\}}$ so, in particular, contains an element $g:=un_8$. Here $u\in U$ and $n_8$ is an element of $N$ which when mapped to the Weyl group is the reflection in root class 8. Without loss of generality we can assume that $g=x_1(t)n_8$. Then
$$gX_1 g^{-1}= x_1(t)n_8X_1 n_8^{-1}x_1(t)^{-1}=x_1(t)X_7x_1(t)^{-1}<U_1.$$
Since $U_1<\GL$ this implies that $X_1<\GL$ which is a contradiction.

\begin{remark}
We are left with the exceptional cases from Lemma \ref{lemma:saxl}. Thus from now on $G$ is a Chevalley group and note that Lemma \ref{lemma:parabolicadj} still applies. In what follows we number the roots in the normal way and refer to parabolic subgroups by the number of the roots which are not included in their generating set. 
\end{remark}

\subsection{$G=PSL_n(q)$}

If $\Galph=P_1$ or $P_{n-1}$ then the action on points is $2$-transitive, $G$ is flag-transitive in its action on $\spaceS$ and the action is well understood. Thus we exclude this possibility and observe that we may assume that $n\geq 4$.

Consider $G$ in the standard projective modular representation. Let $\Galph = P_k$, $k\in\{2,\dots, n-2\}$. By Lemma \ref{lemma:parabolicadj},
$$U_kL_{k-1, k, k+1}H\leq \GL\leq P_{k-1,k+1}.$$

Now without loss of generality $2k\leq n$ (reorder the roots if necessary); then conjugate $\Galph$ by a permutation matrix $g\in G$ corresponding to the $(1,k+1)(2,k+2)\dots(k, 2k)$ permutation.

Then $g\not\in\Galph$ hence $\Galph\cap\Galph^g\leq\GL$. If $n=2k$ this means that $\hatt SL_k(q)\times SL_k(q)\leq \GL$ which is impossible since $\GL\leq P_{k-1, k+1}$. If $n>2k$ then this means that $SL_k(q)\times Q_k\leq \GL$ where $Q_k$ is isomorphic to a $k$-th parabolic group in $SL_{n-k}(q)$. If $k\geq 3$ then this is clearly impossible.

Assume then that $k=2$. We must have $\hatt SL_2(q)\times SL_2(q)\times SL_{n-4}(q)\leq \GL\leq \hatt A:((q-1)\times SL_2(q)\times SL_{n-3}(q))$. Thus either $n=5$ or $SL_2(q)$ is not quasi-simple, i.e. $q=2$ or $3$.

Consider the case when $n=5$. Then 
$$v=(q^2+1)(q^4+q^3+q^2+q+1), \ v-1= q(q^2+q+1)(q^3+q+1).$$
Furthermore $b$ is divisible by $q|G:P_{1,3}| = q(q^2+1)(q^4+q^3+q^2+q+1)(q^2+q+1)$.

Thus $|P_{1,3}:\GL|$ divides into $q(q^3+q+1)$ and so divides into $q(q-1,3)$. In fact $|P_{1,3}:\GL|$ is also divisible by $q$ and $\GL>U_2$. No such subgroup exists for $q>7$. When $q\leq 7$ we must have $k(k-1)$ dividing into $q^3+q+1$. Examining the numerical values of $q^3+q+1$ for $q=2, 3, 5$ and $7$ we find that this is not possible.

We are left with the possibility that $k=2, n\geq 6$ and $q=2$ or $3$. If $q=2$ then conditions on $\GL$ imply that $S_3\times SL_{n-4}(2)<SL_{n-3}(2)$. If $q=3$ we have that $SL_2(3)\times SL_{n-4}(3)<SL_{n-3}(3)\times 2$. In both cases this gives a contradiction.

\subsection{$G=D_m(q)$, $m\geq3$ odd}

If $m=3$ then $G=PSL_4(q)$ and we are already done.

Suppose $m\geq 5$. If $\Galph = P_i, i< m-1$ then Lemma \ref{lemma:saxl} still applies (c.f. \cite[Section 5]{saxl}). The cases where $\Galph=P_m$ or $\Galph=P_{m-1}$ are analogous, so we just consider $\Galph=P_m$. Thus
$$v=(q^{m-1}+1)(q^{m-2}+1)\dots(q^2+1)(q+1).$$
Then Lemma \ref{lemma:parabolicadj} implies that $U_m:L_{m,m-1}\leq\GL\leq P_{m-2}$. Thus $b$ is divisible by
$$(q^{m-1}+\dots+q+1)(q^{m-3}+\dots+q^2+1)\frac{v}{q+1}.$$
If $m\equiv 1(4)$ then $(q^{m-3}+\dots+q^2+1,v)\geq q^2+1.$ If $m\equiv 3(4)$ then $(q^{m-3}+\dots+q^2+1,v)\geq q^{\frac{m-3}{2}}-q^{\frac{m-5}{2}}+\dots-q+1$. When $(m,q)\neq (7,2)$ this contradicts the fact that $b$ divides into $v(v-1)$. A simple combinatoral argument rules out the case when $(m,q)=(7,2)$.

\subsection{$G=E_6(q)$}

If $\Galph = P_i, i=2$ or $4$ then Lemma \ref{lemma:saxl} still applies (c.f \cite[Section 8]{saxl}).

If $\Galph = P_1$ then,
by Lemma \ref{lemma:parabolicadj}, $U_1:L_{1,3}\leq\GL\leq P_3$. This implies that $q^2+1$ divides into $b$. However $(q^2+1, v(v-1))$ divides into $2$. This yields a contradiction.

If $\Galph=P_3$ then, by Lemma \ref{lemma:parabolicadj}, $U_3:L_{1,3,4}\leq\GL\leq P_{1,4}$. This implies that $(q^2+1)^2$ divides into $b$. Now $(v-1, q^2+1)=1$ and $(v/(q^2+1), (q^2+1))\leq 2$. Once again we have a contradiction.

\section{Concluding remarks}

Theorem \ref{theorem:main} has the following corollary: 
\begin{corollary}
Suppose that $G$ has socle $T$ a simple group of Lie type and $G$ acts line-transitively on a linear space $\spaceS$. If a point-stabilizer in $T$ is a parabolic subgroup of $T$ then a line-stabilizer in $T$ is also a parabolic subgroup of $T$.
\end{corollary}

For particular families of low rank simple groups of Lie type Theorem \ref{theorem:main} is implied by existing results in the literature. We have already mentioned the case where $G$ has Lie rank 1; in addition results exist covering the case when $G$ has socle $PSL_3(q)$ \cite{gill2}.

\bibliographystyle{amsalpha}
\bibliography{characteristic}

\providecommand{\bysame}{\leavevmode\hbox to3em{\hrulefill}\thinspace}
\providecommand{\MR}{\relax\ifhmode\unskip\space\fi MR }
\providecommand{\MRhref}[2]{%
  \href{http://www.ams.org/mathscinet-getitem?mr=#1}{#2}
}
\providecommand{\href}[2]{#2}
\begin{thebibliography}{GMS03}

\bibitem[BDD88]{bdd}
F.~Buekenhout, A.~Delandtsheer, and J.~Doyen, \emph{Finite linear spaces with
  flag-transitive groups}, J. Combin. Theory, Series A \textbf{49} (1988),
  268--293.

\bibitem[Car89]{carter}
Roger Carter, \emph{Simple groups of {L}ie type}, John Wiley and Sons, 1989.

\bibitem[CKS76]{csk}
Charles~W. Curtis, William~M. Kantor, and Gary~M. Seitz, \emph{The
  {$2$}-transitive permutation representations of the finite {C}hevalley
  groups}, Trans. Amer. Math. Soc. \textbf{218} (1976), 1--59.

\bibitem[CNP03]{cnp}
A.~Camina, P.~Neumann, and C.~Praeger, \emph{Alternating groups acting on
  linear spaces}, Proc. London Math. Soc.(3) \textbf{87} (2003), no.~1, 29--53.

\bibitem[Gil]{gill2}
Nick Gill, \emph{${PSL}(3,q)$ and line-transitive linear spaces}, Submitted.

\bibitem[GLS94]{gorenstein}
Daniel Gorenstein, Richard Lyons, and Ronald Solomon, \emph{The classification
  of the finite simple groups, number 3}, Mathematical Surveys and Monographs,
  vol.~40, American Mathematical Society, 1994.

\bibitem[GMS03]{gms}
Robert~M. Guralnick, Peter M{\"{u}}ller, and Jan Saxl, \emph{The rational
  function analogue of a question of {S}chur and exceptionality of permutation
  representations}, Mem. Amer. Math. Soc. \textbf{162} (2003), no.~773, 1--79.

\bibitem[KL90]{kl}
P.~Kleidman and M.~Liebeck, \emph{The subgroup structure of the finite simple
  groups}, London Mathematical Society Lecture Note Series, vol. 129, Cambridge
  University Press, Cambridge, 1990.

\bibitem[Kle90]{kleidman4}
Peter~B. Kleidman, \emph{The finite flag-transitive linear spaces with an
  exceptional automorphism group}, Finite geometries and combinatorial designs
  (Lincoln, NE, 1987), Contemp. Math., vol. 111, Amer. Math. Soc., Providence,
  RI, 1990, pp.~117--136.

\bibitem[Mal91]{malle}
Gunter Malle, \emph{The maximal subgroups of ${^2F_4(q^2)}$}, J. Algebra
  \textbf{139} (1991), 52--69.

\bibitem[Sax02]{saxl}
Jan Saxl, \emph{On finite linear spaces with almost simple flag-transitive
  automorphism groups}, J. Combin. Theory A \textbf{100} (2002), 322--348.

\end{thebibliography}

\end{document}